\documentclass[11pt]{amsart}

\usepackage{epsfig,amsmath,amsfonts,latexsym,caption,subcaption,xcolor}
\usepackage[abs]{overpic}
\usepackage{color}
 \usepackage{palatino}
 \usepackage{hyperref}
 \usepackage[nameinlink,noabbrev,capitalise]{cleveref}

\newtheorem{theorem}{Theorem}
\newtheorem{proposition}[theorem]{Proposition}
\newtheorem{lemma}[theorem]{Lemma}

\newtheorem{conjecture}[theorem]{Conjecture}
\newtheorem{question}[theorem]{Question}
\newtheorem*{thmMain}{\cref{main}}
\newtheorem*{thmMain2}{\cref{main2}}
\newtheorem*{thmMain3}{\cref{main3}}
\newtheorem*{thmMain5}{\cref{main5}}
\newtheorem*{lem1}{\cref{simplecobord}}

\theoremstyle{definition}
\newtheorem{definition}[theorem]{Definition}

\newtheorem{example}[theorem]{Example}

\theoremstyle{remark}
\newtheorem{remark}[theorem]{Remark}

\def\dfn#1{{\em #1}}
\def\R{\mathbb{R}}
\def\Z{\mathbb{Z}}

\def\tb{Thurston--Bennequin }
\DeclareMathOperator{\tbn}{tb}
\DeclareMathOperator{\tbm}{\overline{tb}}
\def\co{\colon\thinspace}

\numberwithin{theorem}{section}
\theoremstyle{plain}

\begin{document}

\title{Lagrangian Realizations of Ribbon Cobordisms} 

\author{John B. Etnyre}
\address{School of Mathematics \\ Georgia Institute
of Technology \\  Atlanta  \\ Georgia}
\email{etnyre@math.gatech.edu}

\author{Caitlin Leverson}
\address{Mathematics Program \\ Bard College \\ Annandale-on-Hudson \\ New York}
\email{cleverson@bard.edu}


\begin{abstract}
We show that up to stabilizations, smooth ribbon cobordisms can be realized by decomposable Lagrangian cobordisms. We also define the notion of stabilization for Lagrangian cobordisms and show that it can be used to find new Lagrangian cobordisms between Legendrian links. 
\end{abstract}

\maketitle

\section{Introduction}

Understanding Lagrangian cobordisms between Legendrian knots is a fundamental problem in symplectic and contact topology. There has been a great deal of work constructing such cobordisms \cite{BST,Chantraine10,CornwellNgSivek16,EkholmHondaKalman16,eliashbergGromov}, obstructing them \cite{BLW,BaldwinSivek,capovilla2022,Chantraine10,Chantraine12,chantraine2020floer,CornwellNgSivek16,EkholmHondaKalman16,GollaJuhasz,Pan17augCat,SabloffTraynor}, and distinguishing them \cite{CasalsGao22,CasalsNg,Pan17torus}. They are also a central part of symplectic field theory as they induce maps between the contact homology of Legendrian knots \cite{EkholmHondaKalman16,Pan17augCat}. In this article we will give several general existence results for Lagrangian cobordisms up to stabilization, analogous to the well known fact that every smooth knot in $S^3$ has a Legendrian realization (in fact, infinitely many).

\subsection{Context}
We begin by recalling that if $(Y_\pm,\xi_\pm)$ are closed contact $3$-manifolds, then a \dfn{symplectic cobordism} from $(Y_-,\xi_-)$ to $(Y_+,\xi_+)$ is a compact symplectic $4$-manifold $(X,\omega)$ with $\partial X=-Y_-\cup Y_+$, and $(Y_-,\xi_-)$ is a concave boundary component and $(Y_+,\xi_+)$ is a convex boundary component. We will frequently denote $Y_\pm$ by $\partial_\pm X$. If $\Lambda_\pm$ is a Legendrian link in $(Y_\pm,\xi_\pm)$ then a \dfn{Lagrangian cobordism} from $\Lambda_-$ to $\Lambda_+$ is a Lagrangian surface $L$ in $X$ that is transverse to $\partial X$ and $\partial L=-\Lambda_-\cup \Lambda_+$. We note that most definitions in the literature require that $L$ be a product near the boundary \cite{Chantraine10, EkholmHondaKalman16}, but our first result (originally worked out in joint work with the first author, Lenhard Ng, and Tye Lidman, but never published) shows that this is not necessary and is implied by the fact that $\partial L$ is Legendrian and that $L$ is transverse to $\partial X$. 
\begin{lemma}\label{simplecobord}\label{maketangent}
Given a Lagrangian cobordism $L$ in a strong symplectic cobordism $(X,\omega)$ there are Liouville vector fields $v_i$ defined near $\partial X$ that are tangent to $L$. 
\end{lemma}
If $(X,\omega)$ is an exact symplectic manifold, so that $\omega=d\beta$ for some $1$-form $\beta$, then we say that $L$ is an \dfn{exact Lagrangian cobordism} if in addition to the requirements above, $\beta|_L=df$ for some function $f$ with the property that $f$ is constant on $\Lambda_-$ and is constant on $\Lambda_+$. One may easily check that $f$ must be constant on each component of $\Lambda_\pm$ and thus the latter part of the definition is automatic if $\Lambda_-$ and $\Lambda_+$ are knots. However, this condition is important when considering multi-component links \cite{Chantraine15note}. We will be interested in exact Lagrangian cobordisms in the symplectization of the standard contact $S^3$ which is the symplectic manifold $([0,1]_t\times S^3,\omega)$ where $\omega=d(e^t\alpha)$ for $\alpha$ the standard contact form on $S^3$ and where $t$ is the coordinate on $[0,1]$.

\subsection{Main Results}
One basic question about exact Lagrangian cobordisms is: Given two Legendrians $\Lambda_\pm$, when does there exist an exact Lagrangian cobordism from $\Lambda_-$ to $\Lambda_+$? There are a number of known obstructions coming from, for example, the classical Legendrian knot invariants \cite{Chantraine10}, contact homology of Legendrian knots \cite{capovilla2022,EkholmHondaKalman16,Pan17augCat}, ruling polynomials \cite{Pan17augCat}, generating families \cite{SabloffTraynor}, and Heegaard Floer Theory \cite{BLW,BaldwinSivek,GollaJuhasz}. (See \cite{Blackwell2021} for a survey of obstructions to the existence of Lagrangian cobordisms.) There are also a number of known constructions. In particular, it is known that if $\Lambda_-$ and $\Lambda_+$ are Legendrian isotopic then there exists a Lagrangian cobordism from $\Lambda_-$ to $\Lambda_+$ \cite{Chantraine10,EkholmHondaKalman16,eliashbergGromov}. From Ekholm, Honda, and K\'{a}lm\'{a}n \cite{EkholmHondaKalman16}, we also know there exists a Lagrangian cobordism from $\Lambda_-$ to $\Lambda_+$ if $\Lambda_-$ can be obtained from $\Lambda_+$ by a ``pinch" move (see \cref{band}) or if $\Lambda_+=\Lambda_-\cup U$ where $U$ is a Legendrian unknot with maximal Thurston-Bennequin number which can be ``separated" from $\Lambda_-$ by a topological ball (see \cref{birth}). Note that smoothly these correspond to isotopy, saddles, and minima, respectively. Cobordisms which are the result of stacking such pieces are called \dfn{decomposable Lagrangian cobordisms}. See \cref{dlc} for the precise definition. As none of the above cobordisms are maxima when considered smoothly, all decomposable Lagrangian cobordisms are ribbon cobordisms (see \cref{ribbon} for the definition of a ribbon cobordism). A natural question is then, given a ribbon cobordism $C$, does there exist a decomposable Lagrangian cobordism which is isotopic to $C$? This article answers the question in the affirmative. 

\begin{theorem}\label{main}
Given a ribbon cobordism $C$ in $[0,1]\times S^3$ from the link $K_-$ to the link $K_+$ with $K_-\neq\emptyset$ and a Legendrian realization $\Lambda_\pm$ of $K_\pm$, then $C$ may be isotoped to a decomposable Lagrangian cobordism $L$ from some stabilization of $\Lambda_-$ to some stabilization of $\Lambda_+$. 
\end{theorem}

\begin{remark}
Note that because $C$ is a ribbon cobordism, the assumption that $K_-\neq\emptyset$ implies we also have $K_+\neq\emptyset$.
\end{remark}
\begin{remark}
In \cite{Rizell2024pre}, Dimitroglou Rizell explores similar questions for concordances. In particular, given a totally real concordance $C$, after stabilizing the Legendrians sufficiently many times, one can isotope $C$ to a Lagrangian concordance.
\end{remark}

\cref{main} should be interpreted as part of the general understanding that when allowing ``stabilizations", the rigidity of symplectic and contact geometry becomes flexible and controlled by the topology of the situation. For example, two Legendrian knots are Legendrian isotopic after some number of stabilizations if and only if they are smoothly isotopic, so the subtle distinctions between different Legendrian representatives of a fixed knot type vanish after sufficiently many stabilizations.  

While the previous theorem tells us we can Lagrangian realize every ribbon cobordism, it does not allow us to specify the Legendrian knot types of $\Lambda_\pm$. Our second theorem says that we do not need to stabilize $\Lambda_+$ (and thus can fix the Legendrian knot type of $\Lambda_+$) if we are willing to change the smooth type of the cobordism.

\begin{theorem}\label{main4}\label{main2}
Given a ribbon cobordism $C$ in $[0,1]\times S^3$ from the link $K_-$ to the link $K_+$ with $K_-\neq\emptyset$ and a Legendrian realization $\Lambda_+$ of $K_+$, there exists a Legendrian realization $\Lambda_-$ of $K_-$ and a decomposable Lagrangian cobordism $L$ from $\Lambda_-$ to $\Lambda_+$. (The topology of $L$ might be different from that of $C$.)
\end{theorem}

The proof of this theorem relies on a notion of stabilization for a Lagrangian cobordism $L$ from $\Lambda_-$ to $\Lambda_+$, which we define in \cref{stabilize}. More specifically, given an arc $\gamma$ properly embedded in $L$ with $\gamma\cap \Lambda_\pm$ a single point in each of $\Lambda_-$ and $\Lambda_+$, we can define the $\pm$-stabilization $S_\pm^\gamma(L)$. This stabilization of $L$ will be a Lagrangian cobordism from $S_\pm(\Lambda_-)$ to $S_\pm(\Lambda_+)$ and will only depend on $\gamma$ up to isotopy in $L$. This form of stabilization also makes the Lagrangian cobordisms ``flexible".

 We believe the following more general conjecture will follow from techniques similar to those used in our proofs.
\begin{conjecture}\label{maincob}
If $L_1$ and $L_2$ are two decomposable Lagrangian cobordisms that are smoothly ribbon isotopic, then after stabilizing both cobordisms, they become Lagrangian isotopic. 
\end{conjecture}

\begin{remark}
We note, for simplicity of the statements and the proofs, that the theorems above have been stated for the symplectic manifold $[0,1]\times S^3$, but it is not much harder to prove the same results for Stein cobordisms. We do not know if they are true in general symplectic cobordisms. 
\end{remark}

We also ask the following question about simplifying cobordisms by stabilizing them. 
\begin{question}
Can every Lagrangian ribbon cobordism from $L_1\not=\emptyset$ to $L_2\not=\emptyset$ be stabilized to become decomposable?
\end{question}

Analogous to Lagrangian cobordisms between Legendrian knots, we can also discuss symplectic cobordisms between transverse knots. In this context, we have analogs of Theorems~\ref{main} and~\ref{main2}. Specifically, the analog of \cref{main} is the following. 

\begin{theorem}\label{main3}
Given a ribbon cobordism $C$ in $[0,1]\times S^3$ from the link $K_-$ to the link $K_+$ with $K_-\neq\emptyset$, then there is some transverse realization $T_\pm$ of $K_\pm$ so that $C$ may be isotoped to a symplectic cobordism from $T_-$ to $T_+$. 
\end{theorem}

The analog of \cref{main2} is the following.
\begin{theorem}\label{main5}
Given a ribbon cobordism $C$ in $[0,1]\times S^3$ from the link $K_-$ to the link $K_+$ with $K_-\neq\emptyset$ and a transverse realization $T_+$ of $K_+$, then there is some transverse realization $T_-$ of $K_-$ so that $C$ may be isotoped to a symplectic cobordism from $T_-$ to $T_+$. 
\end{theorem}
We are not sure if \cref{maincob} has an analog in the symplectic case.

\subsection*{Organization} In \cref{stabcobord} we define the notion of a stabilization of a Lagrangian cobordism. \cref{ribbon} recalls the definition of a ribbon cobordism as well as defining a normal form for a ribbon cobordism which is used in the proof of \cref{main}. We recall the definition of a decomposable Lagrangian cobordism and prove Theorems~\ref{main} and~\ref{main3} in \cref{dlc}. In \cref{neighbor} we prove \cref{simplecobord} and Theorems~\ref{main3} and~\ref{main5}. Finally, in \cref{examples} we explore some examples.

\subsection*{Acknowledgements} The authors thank Lenhard Ng and Tye Lidman for allowing them to include the proof of \cref{simplecobord} in this manuscript and also thank Lenhard Ng for helpful conversations. The first author was partially supported by National Science Foundation grants DMS-1906414 and DMS-2203312 and the Georgia Institute of Technology Elaine M. Hubbard Distinguished Faculty Award. The second author was supported by NSF Postdoctoral Fellowship DMS-1703356.

\section{Stabilizing Lagrangian cobordisms}\label{stabcobord}

In this section we will define the stabilization of a Lagrangian cobordism. We begin by observing that a properly embedded arc in a Lagrangian cobordism has a standard neighborhood that looks like the symplectization of a neighborhood of a point on a Legendrian knot.

\begin{lemma}\label{stab1}
Let $L$ be a Lagrangian cobordism in $(X,\omega)$ from $\Lambda_-\neq\emptyset$ to $\Lambda_+\neq\emptyset$ and let $\gamma:[0,1]\to L$ be a properly embedded arc in $L$ connecting a point on $\Lambda_-$ to a point on $\Lambda_+$. Then there is a neighborhood $N$ of the image of $\gamma$ in $X$ that is symplectomorphic to the symplectization $[0,1]\times N'$ of a neighborhood $N'$ of the Legendrian arc $\Lambda=[-\epsilon, \epsilon]$ in $\partial_-X$ so that  $L\cap (\partial_\pm X)\cap N$ is mapped to $\Lambda_\pm$ and $L\cap N$ is taken to a portion of $[0,1]\times \Lambda$.
\end{lemma}
\begin{proof}
A slight variant on the standard Moser-type argument, see \cite[Proposition~4.2]{Weinstein91}, shows that there is a contactomorphism, that preserves a given contact form, from a neighborhood $N_-$ of $\gamma(0)$ in $\partial_- X$ and a neighborhood $N_+$ of $\gamma(1)$ in $\partial_+X$  to a neighborhood of $(0,0,0)$ in $\{i\}\times \R^3$ that takes $\Lambda_\pm\cap N_\pm$ to a portion of the $x$-axis that we denote by $\Lambda$.  Using the Liouville vector fields for the symplectizations, these neighborhoods can be extended to neighborhoods $N'_\pm$ in $X$ that are symplectomorphic to neighborhoods in $[0,1]\times \R^3$ by a symplectomorphism that takes $L\cap N'_\pm$ to a portion of the $tx$-plane.  We now extend these symplectomorphisms to a diffeomorphism from a neighborhood $N'$ of $\gamma([0,1])$ to a neighborhood of $[0,1]\times\{(0,0,0)\}$ that takes $N\cap L$ to $[0,1]\times\Lambda$. Since $L$ and $[0,1]\times \Lambda$ are Lagrangian, another standard Moser-type argument will isotope the diffeomorphisms into a symplectomorphism on a possibly smaller neighborhood $N$. This neighborhood and symplectomoprhism satisfy all the claimed conditions. 
\end{proof}

This allows us to define the notion of a stabilization of a Lagrangian cobordism.

\begin{definition}[Stabilization of a Lagrangian cobordism]\label{stabilize}
Given a Lagrangian cobordism $L$ in $X$ from $\Lambda_-\neq\emptyset$ to $\Lambda_+\neq\emptyset$ and a properly embedded arc $\gamma:[0,1]\to L$ in $L$ connecting a point on $\Lambda_-$ to a point on $\Lambda_+$, we define the \dfn{$\pm$-stabilization of $L$ along $\gamma$}, denoted $S^\gamma_\pm(L)$, as follows: let $N$ be the neighborhood of the image of $\gamma$ from the previous lemma and replace $L\cap N$ by $[0,1]\times S_\pm(\Lambda)$.
\end{definition}

The following implies that $S^\gamma_\pm(L)$ does not depend on any choices when $L$ is a concordance. 

\begin{lemma}
The $\pm$-stabilization of a Lagrangian cobordism $L$ along $\gamma$ only depends on the isotopy class of $\gamma$.
\end{lemma}

\begin{proof}
The Moser-type arguments in \cref{stab1} can be done in $1$-parameter families. So, given an isotopy $\gamma_t$ we have a family of neighborhoods $N_t$ and the stabilization construction will give a Legendrian isotopy. 
\end{proof}

\section{Ribbon cobordisms}\label{ribbon}
Recall that a properly embedded surface $C$ in $[0,1]\times S^3$ with $\partial C= -K_0 \cup K_1$ where $K_i=C\cap (\{i\}\times S^3)$ is called a \dfn{ribbon cobordism} if the projection map $p:[0,1]\times S^3\to [0,1]$ restricted to $C$ is Morse and only has index $0$ and $1$ critical points. If the genus of $C$ is zero, then we call it a \dfn{ribbon concordance}. 

If $K_0=\emptyset$ we call $C$ a \dfn{ribbon surface} for $K_1$ and notice that it is equivalent to a surface immersed in $S^3$ with only ribbon singularities.

\begin{remark}
We note that if $X$ is a Stein manifold with plurisubharmonic function $\phi\co X\to \R$, then we may make the same definition with $\phi$ in place of the projection to $[0,1]$ above. For simplicity, we will focus our discussion on $[0,1]\times S^3$.
\end{remark}

\begin{figure}[htbp]\tiny
    \begin{center}
    \begin{overpic}
    {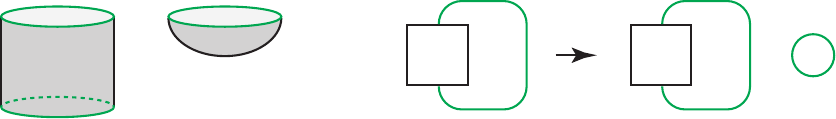}
      \put(204, 26){\Large$K$}
      \put(311, 26){\Large$K$}
    \end{overpic}
    \caption{The birth move. On the left the cobordism near an index $0$ critical point ($t$ is increasing as you move up). On the right is the effect on the knot diagram when one passes an index $0$ critical point in a cobordism as $t$ increases in the direction of the arrow. }
    \label{birth}
    \end{center}
  \end{figure}

\begin{figure}[htbp]\small
    \begin{center}
    \begin{overpic}
    {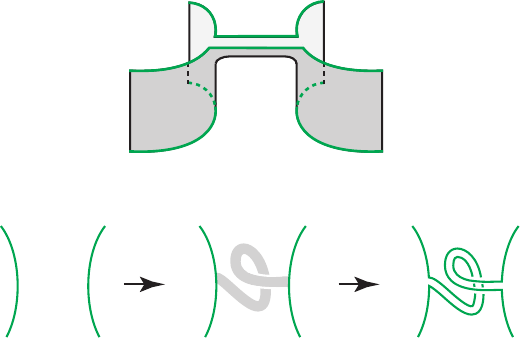}
      \put(11, 59){$t_p-2\epsilon$}
      \put(111, 59){$t_p-\epsilon$}
      \put(221, 59){$t_p$}
	\put(11, 5){$L_p^-$}
	\put(250, 5){$L_p^+$}
    \end{overpic}
    \caption{The band move. On the top is a schematic of the cobordism near an index $1$ critical point ($t$ is increasing as you move up). On the bottom is the effect on the knot diagram when one passes an index $1$ critical point in a cobordism. The middle figure is the flattened band about the unstable manifold $\gamma$. This band can twist (specified by a framing) and knot.}
    \label{band}
    \end{center}
  \end{figure}
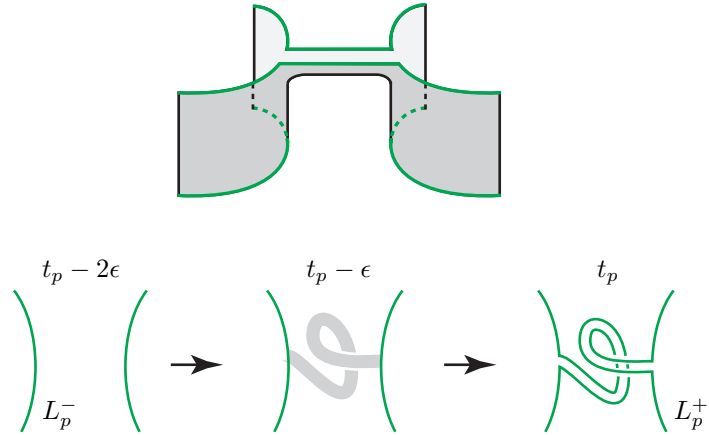

Given a ribbon cobordism we can isotope it into a normal form as follows. First, we can arrange for all the index $0$ critical points to occur before all the index $1$ critical points, and for all the critical points to occur at distinct $t$ values (where $t$ is the coordinate on $[0,1]$). So $C$ can be constructed by moving from $\{0\}\times S^3$ to $\{1\}\times S^3$ by applying a sequence of the following moves, see Figures~\ref{birth} and~\ref{band}:
\begin{enumerate}
\item Traces of isotopies of links in $S^3$.
\item Births of unknots (disjoint from the rest of the diagram).
\item Band moves. 
\end{enumerate}
The births of unknots occur when an index $0$ critical point is passed. A band move occurs when an index $1$ critical point is passed. To better understand the band move, consider an index $1$ critical point $p$ that occurs at time $t_p$. Let $\gamma_p$ be the unstable manifold of $p$ in $C\cap ([t_p-\epsilon, t_p]\times S^3)$, where $\epsilon$ is chosen small enough that there are no other critical values in $[t_p-\epsilon, t_p]$. We can ``flatten'' the $1$--handle corresponding to $p$ by pushing a small neighborhood of $\gamma_p$ to lie entirely in $\{t_p-\epsilon\}\times S^3$. Now just before $t_p-\epsilon$ we see a link $L_p^-$ in $S^3$, at $t_p-\epsilon$ a band is attached to this link, and just after $t_p-\epsilon$ we have the result $L_p^+$ of the band move. See \cref{band}. Specifically, let $B$ be the neighborhood of $\gamma$ in $C$ and set $\partial_- B=B\cap L_p^+$ and $\partial_+ B= \overline{(\partial B-\partial_-B)}$.

Notice that the band moves are determined entirely by the curves $\gamma_p$ and a relative framing to indicate how $\gamma$ is thickened into a band, and if $\gamma_p$ is isotoped through curves with end points on $L_p^-$ and interiors disjoint from $L_p^-$ then the resulting surface is the same, up to isotopy. It is possible that a curve $\gamma_q$ corresponding to the critical point $q$ could have its end points on part $L_q^-$ corresponding to a band attached earlier. But we can always isotop $\gamma_q$ so that its endpoints are on part of $L_q^-$ corresponding to $K_0$ and the boundaries of the birth disks. 

The following lemma is probably well-known to experts, but we could not find a reference and so include it here for the convenience of the reader. 
\begin{lemma}\label{reorder}
If the ribbon cobordism $C$ is built as above and $\gamma_p$ and $\gamma_q$ correspond to adjacent band moves with $\gamma_p$ being attached first, then $C$ can be isotoped so that the band $\gamma_q$ is attached before $\gamma_p$ and nothing else about the presentation of the cobordism changes. 
\end{lemma}
\begin{remark}
The discussion above together with this lemma says every ribbon cobordism is isotopic to one built by some number of births of unknots, followed by a sequence of band moves and the bands can be attached in any oder. (And between any of these moves, we can also look at the traces of isotopies.)
\end{remark}
\begin{proof}
Recall from above, that we are assuming that the $\gamma_q$ has its end points disjoint from the parts of $L_q^+$ that correspond to previous bands. 
Let $B_p$ and $B_q$ be the bands corresponding to the band that are neighborhoods of $\gamma_p$ and $\gamma_q$, respectively and recall that we can take them as thin as we like. So if, when we push $\gamma_q$ to the level that $B_p$ is attached, it is disjoint from $B_p$, then we can also assume $B_q$ is and further push it down to be attached before $B_p$. 

So the only concern is when $\gamma_q$ intersects $B_p$. In this case consider $L_p^+=L_q^-$ just after $B_p$ was attached. Recalling where $B_p$ was, we can push any intersection point of $\gamma_q$ off of $B_p$ through $\partial_-B_p$. See \cref{fix}. We can now push $\gamma_q$ down past the level of $B_p$.  
  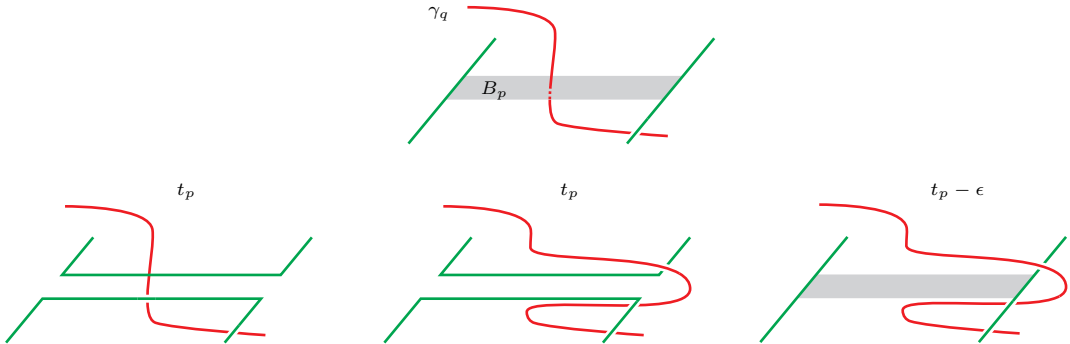
\begin{figure}[htbp]
    \begin{center}\tiny
    \begin{overpic}
    {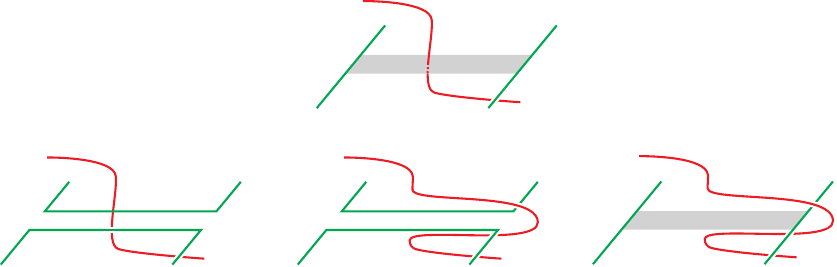}
    \put(65, 57){$t_p$}
    \put(210, 57){$t_p$}
    \put(350, 57){$t_p-\epsilon$}
      \put(180, 95){$B_p$}
      \put(160, 125){$\gamma_q$}
    \end{overpic}
    \caption{On the top we see an example obstructing us from pushing the band from $q$ down past the band for $p$. In the next row, we see how to fix this problem. Specifically, at a level above $t_p-\epsilon$ we isotope $\gamma_q$ as discussed in the proof so that when pushed to the level $B_p$ is attached, $\gamma_q$ is disjoint from it.}
    \label{fix}
    \end{center}
  \end{figure}
\end{proof}

Though not essential for our applications, it is simpler to understand our construction of decomposable Lagrangian cobordisms if we assume there are no isotopies involved in a ribbon cobordism. 
\begin{lemma}\label{noiso}
Given a ribbon cobordism $C$ from $K_0$ to $K_1$, after isotoping $K_1$ and $C$ we can assume that $C$ is built from product cobordisms, births, and band moves. 
\end{lemma}
\begin{proof}
Suppose that there are no critical values between $t_1$ and $t_2$, with $t_1<t_2$. Then by definition $C$ corresponds to the trace of an isotopy of a link $K$. That is, there are links $K_t$ for $t\in[t_1,t_2]$ such that $C$ is the union of $(t,K_t)$. Now the Isotopy Extension Theorem says there is an isotopy $\phi_t$ of $M$ for $t\in [t_1,t_2]$, such that $K_t=\phi_{t-t_1}(K_1)$. We can now apply $\phi_t$ to $\{t\}\times M$ and $\phi_{t_n}$ to all the $\{t\}\times M$ with $t>t_2$. This will result in a cobordism that is constant in the region $[t_1,t_2]\times M$. 

Thus outside a small neighborhood of the critical levels, the cobordism can be assumed to be constant (though of course $K_1$ will have been isotoped from its initial position). Now in a neighborhood of a critical point, we can apply the same argument outside the neighborhood of the critical point to show that the cobordism is constant except in the neighborhood where the birth or band move is performed. 
\end{proof}

\section{Decomposable Lagrangian cobordisms}\label{dlc}
A \dfn{decomposable Lagrangian cobordism} from a Legendrian link $\Lambda_-$ to a Legendrian link $\Lambda_+$ in a contact manifold $(M,\xi)$ is a Lagrangian submanifold of a piece of the symplectization of $(M,\xi)$ that is made from a sequence of the following moves. See \cref{decomp}.
\begin{enumerate}
\item Traces of a Legendrian isotopy\footnote{Actually one must perturb the traces to obtain a Lagrangian concordance.}.
\item Births of a maximal Thurston-Bennequin unknot (disjoint from the rest of the diagram).
\item Pinch moves. 
\end{enumerate}
  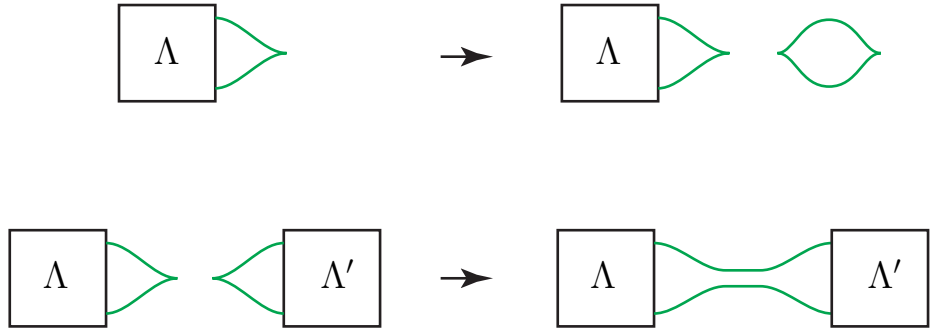
\begin{figure}[htbp]\tiny
    \begin{center}
    \begin{overpic}
    {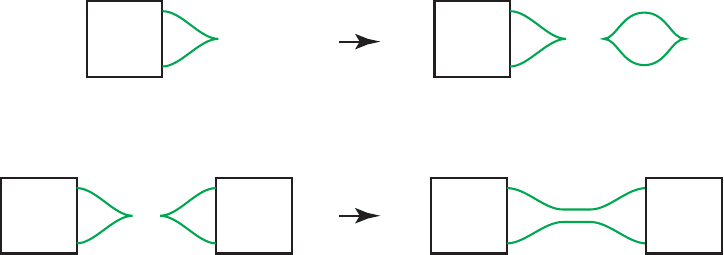}
      \put(13, 15){\Large$\Lambda$}
      \put(118, 15){\Large$\Lambda'$}
      \put(220, 15){\Large$\Lambda$}
      \put(325, 15){\Large$\Lambda'$}
      \put(55, 100){\Large$\Lambda$}
      \put(222, 100){\Large$\Lambda$}
    \end{overpic}
    \caption{In the top row is the birth of a maximal Thurston-Bennequin unknot in a cobordism. On the bottom row is a pinch move (here $\Lambda$ and $\Lambda'$ can be the same component or link). All arrows indicate increasing $t$-values.}
    \label{decomp}
    \end{center}
  \end{figure}

\begin{figure}[htbp]\tiny
    \begin{center}
    \begin{overpic}
    {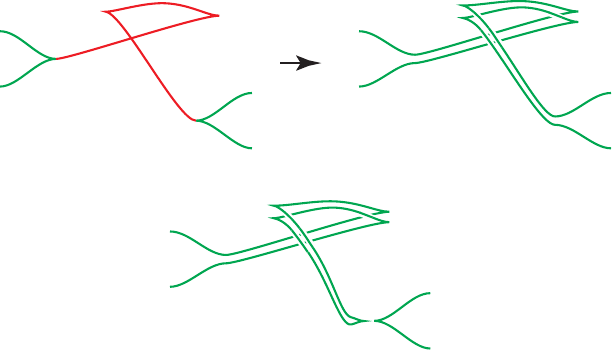}
     \put(76, 140){$\gamma$}
    \end{overpic}
    \caption{In the top row, a more general pinch move in front projections that is equivalent to the original one up to Legendrian isotopy. In the bottom row the Legendrian isotopy proving this.}
    \label{gpinch}
    \end{center}
  \end{figure}

We notice that up to Legendrian isotopy we can describe the pinch move in another way. Let $\gamma$ be any Legendrian arc with end points on a Legendrian $\Lambda$. A standard Moser type argument shows that there is a neighborhood $N$ of $\gamma$ such that it is contactomorphic to a neighborhood of the origin in the standard contact structure on $\R^3$ and $(\gamma\cup \Lambda)\cap N$ maps to a neighborhood of the cusp in the bottom left of \cref{decomp} union a horizontal line connecting them. Thus we may perform a pinch move along $\gamma$. If $\Lambda$ and $\gamma$ are in $(\R^3,\xi_{std})$ and the endpoints of $\gamma$ are cusps in the front projection of $\Lambda$, then the pinch move is realized by taking a copy of $\gamma$ pushed slightly up in the $z$-direction and another pushed slightly down, and connecting them to $\Lambda$ with a small neighborhood of $\gamma$ removed. See \cref{gpinch}.
  
We are now ready to prove our main theorem, which we restate here for convenience.
\begin{thmMain}
Given a ribbon cobordism $C$ in $[0,1]\times S^3$ from the link $K_-$ to the link $K_+$ with $K_-\neq\emptyset$ and a Legendrian realization $\Lambda_\pm$ of $K_\pm$, then $C$ may be isotoped to a decomposable Lagrangian cobordism $L$ from some stabilization of $\Lambda_-$ to some stabilization of $\Lambda_+$. 
\end{thmMain}

\begin{remark}
Note that because $C$ is a ribbon cobordism, the assumption that $K_-\neq\emptyset$ implies we also have $K_+\neq\emptyset$.
\end{remark}

\begin{proof}[Proof of \cref{main}]
We begin by showing that if there is a cobordism from $K$ to $K'$ consisting of a single band move, then if $\Lambda$ is a Legendrian realization of $K$ and one of the components of $\Lambda$ that the band is attached to is stabilized sufficiently many times, then there is a pinch move from $\Lambda$ to a Legendrian realization of $K'$. To see this, let $\gamma$ be the core of the band corresponding to the band move. We can isotope $\gamma$ so that it has end points at one left cusp and one right cusp of $\Lambda$. We can then Legendrian approximate $\gamma$ by $\gamma'$. Now the pinch move along $\gamma'$ corresponds to a band move. It will give a Legendrian realization of $K'$ if the contact framing of $\gamma'$ agrees with the framing of the band move being done.  This does not need to be the case, but can be fixed as follows. If the band move framing is less than the contact framing of $\gamma'$ then we stabilize $\gamma'$ to reduce the contact framing. Each stabilization adds one full left handed twist to the contact framing. So we can arrange that the contact framing agrees with the band move framing or is one less than it (recall the framing on a band can be any half integer). We now are in the situation where the framings agree or the band move framing is larger than the contact framing. If we are in the latter case then near the cusp of the component of $\Lambda$ that has been stabilized many times we can isotope $\Lambda$ as indicated in \cref{fixframing}. 
  \begin{figure}[htbp]\tiny
    \begin{center}
    \begin{overpic}
    {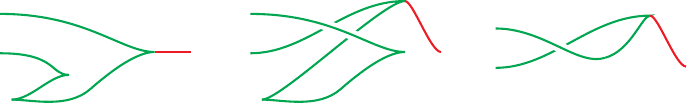}
     \put(50, 4){$\Lambda$}
     \put(80, 29){$\gamma'$}
     \put(172, 4){$\Lambda$}
     \put(206, 36){$\gamma'$}
     \put(300, 19){$\Lambda$}
     \put(325, 33){$\gamma'$}
    \end{overpic}
    \caption{On the left is one end of $\gamma'$. In the middle we have isotoped $\Lambda$ and $\gamma'$. The right most figure is a further isotopy and the framing of a band attached along $\gamma'$ with contact framing will be increased by $1/2$ compared to the band attached to the left most picture.}
    \label{fixframing}
    \end{center}
  \end{figure}
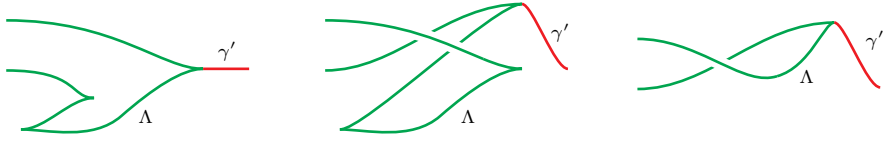
Attaching the result of doing a pinch move along $\gamma'$ now will be the same as doing the same band move before the isotopy, except with framing increased by $1/2$. If there are sufficiently many stabilizations we can increase the framing for the pinch move to agree with that of the band move. 

Now let $C$ be a ribbon cobordism from $K_-$ to $K_+$ in $[0,1]\times S^3$. After isotopy we can assume the cobordism is a sequence of product cobordisms, births, and band moves, by \cref{noiso}; moreover, all the births occur before the band moves. Let $K'$ be the result of all the births in the cobordism $C$ and let $C_1$ be the cobordism from $K_-$ to $K'$. Now let $B_1,\ldots, B_k$ be all of the band moves, labeled in the order that they occur from lowest $t$-coordinate to highest. For each $i$, let $K'_i$ be the knot obtained from $K'$ by attaching the bands $B_1, B_2, B_3,$ through $B_i$. 

Given that the portion of the cobordism $C$ from $K_-$ to $K'$ is purely births, clearly, this cobordism can be made a decomposable cobordism from any Legendrian realization $\Lambda'_-$ of $K_-$ to some Legendrian realization $\Lambda'$ of $K'$. We now want to continue to turn the band moves $B_i$ into pinch moves. To this end, we notice as long as the $B_i$ has one of its endpoints on a component of a Legendrian realization of $K'_{i-1}$ that has been stabilized sufficiently many times then our observation above allows us to realize $B_i$ as a pinch move. By reordering the bands using \cref{reorder} we can assume that each $B_i$ does have an endpoint on a component of a Legendrian realization of $K'_{i-1}$ that was connected to one of the original components of $K_-$. So if the Legendrian realization of $K_-$ has been stabilized sufficiently many times, then all the band moves can be realized by pinch moves and so we have constructed a decomposable Lagrangian cobordism from some Legendrian realization $\Lambda'_-$ of $K_-$ to some realization $\Lambda'_+$ of $K_+$. 

Now if $\Lambda_-$ and $\Lambda_+$ were given Legendrian realizations of $K_-$ and $K_+$ then we notice that some stabilization of $\Lambda_-$ is Legendrian isotopic to $\Lambda'_-$. Thus after stabilizing (as defined in \cref{stabcobord}) the cobordism constructed above sufficiently many times, we can assume that the lower boundary of the cobordism is isotopic to a stabilization of $\Lambda_-$. We can similarly arrange for the top of the cobordism to be isotopic to a stabilization of $\Lambda_+$. 
\end{proof}

We now see in the construction above that we can fix the Legendrian representative of $K_+$ at the expense of stabilizing the cobordism. 
\begin{thmMain2}
Given a ribbon cobordism $C$ in $[0,1]\times S^3$ from the link $K_-$ to the link $K_+$ with $K_-\neq\emptyset$ and a Legendrian realization $\Lambda_+$ of $K_+$, there exists a Legendrian realization $\Lambda_-$ of $K_-$ and a decomposable Lagrangian cobordism $L$ from $\Lambda_-$ to $\Lambda_+$. (The topology of $L$ might be different from that of $C$.)
\end{thmMain2}

\begin{proof}[Proof of \cref{main2}]
In \cite{SabloffVela-VickWong21pre} it was shown how to build a genus 1 Lagrangian cobordism from $S_+(S_-(\Lambda))$ to $\Lambda$. See Figure~\ref{topstab}. 
  \begin{figure}[htbp]\tiny
    \begin{center}
    \begin{overpic}
    {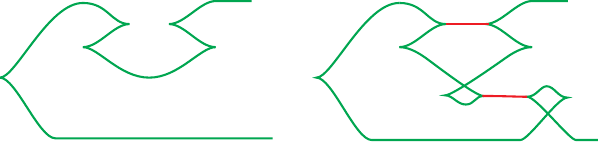}
    \end{overpic}
    \caption{Performing a pinch move along the two red arcs on the right will create a genus 1 Lagrangian cobordism from $S_+(S_-(\Lambda))$ to $\Lambda$.}
    \label{topstab}
    \end{center}
  \end{figure}
Given the Lagrangian cobordism from a stabilization of $\Lambda_-$ to a stabilization of $\Lambda_+$ from \cref{main}, we can further stabilize the cobordism, as defined in \cref{stabcobord}, to obtain a decomposable cobordism with the top being $(S^n_+\circ S^n_-)(\Lambda_+)$ for some $n$. Now applying the cobordism in \cref{topstab} $n$ times, we see a decomposable cobordism from a stabilization of $\Lambda_-$ to exactly $\Lambda_+$ and that this cobordism is smoothly isotopic to a smooth stabilization of $C$. 
\end{proof}

\section{Neighborhoods of Lagrangian cobordisms} \label{neighbor}
In this section we will prove that ribbon cobordisms can be approximated by symplectic cobordisms between transverse knots. To this end, we begin by proving a regular neighborhood theorem for Legendrian submanifolds and show how these Legendrian submanifolds can be perturbed into symplectic manifolds. We first prove \cref{maketangent} which was proven by the first author, Lenhard Ng, and Tye Lidman, but never published.

\begin{lem1}
Given a Lagrangian cobordism $L$ in a strong symplectic cobordism $(X,\omega)$ there are Liouville vector fields $v_i$ defined near $\partial X$ that are tangent to $L$. 
\end{lem1}

\begin{proof}
Let $\partial X=-Y_-\cup Y_+$ and $\Lambda_\pm=L\cap Y_\pm$. We will focus on showing that $v_-$ can be chosen to be tangent to $L$ near $Y_-$, the discussion for $v_+$ is similar. Let $L'$ be the trace of the flow of $\Lambda_-$ under $v_-$. Clearly $L'$ is a Legendrian annulus and $L'$ agrees with $L$ in $Y_-$. Choosing coordinates $(x,\theta)$ on $L'=[0,\epsilon)\times S^1$ there is an identification of a neighborhood of $L'$ in $X$ with a neighborhood $N$ of the zero section in $T^*L'$ that preserves the symplectic structure. We can take coordinate $(x,\theta, y, z)$ on $N=[0,\epsilon)\times S^1\times D^2$ so that the symplectic structure is $\omega= d\lambda$ where $\lambda = y\, dx+ z\, d\theta$. 

Since $L$ is transverse to $Y_-$ we can (near $\Lambda_-$) write $L$ as a graph over $L'$ in $N$. More specifically, we have functions $f,g:[0,\epsilon)\times S^1 \to \R$ such that $L$ near $Y_-$ is parameterized by the map $h(t, \phi)= (x,\theta, f(x,\theta), g(x,\theta))$. Consider 
\[
\beta=h^*\lambda= f\, dx + g\, d\theta.
\]
Since $L$ is Lagrangian we clearly have $d\beta=0$. Moreover, since $L=L'$ when $t=0$ we see that $\beta_{\{(0,\phi)\}}=0$. Thus $\beta$ represents the zero element in the cohomology of $L'$. So there is a function $F:L'\to \R$ so that $dF=\beta$. If $\rho$ is a function on $L'$ that is $0$ near $\Lambda_-$ and $1$ outside a small neighborhood of $\Lambda_-$ then $\beta_s=sdF + (1-s) d(\rho F)$ is a family of closed 1-forms on $L'$ that can be written as $f_s\, dx+ g_s\, d\theta$. The functions $h_s(x,\theta)= (x,\theta, f_s(x,\theta), g_s(x,\theta))$ parameterize a $1$-parameter family of Lagrangian annuli that all agree with $L$ outside a neighborhood of $\Lambda_-$. Thus we get a family of Lagrangian cobordisms $L_s$ such that $L_0=L$ and $L_1$ agrees with $L'$ near $\Lambda_-$. In addition, the $L_s$ are fixed outside a small neighborhood of $\Lambda_-$. 

Choosing an almost complex structure $J$ on $X$ that is compatible with $\omega$ we can apply $J$ to a trivialization of the tangent bundle of $[0,1]\times S^1$ to get a trivialization of the tangent bundles of the $L_s$. We can use the trivialization of the normal bundle to build an isotopy of neighborhoods of the $L_s$ and then extend this to an ambient smooth isotopy $\phi_s:X \to X$ that is fixed on $\partial X$ and outside a neighborhood of $\Lambda_-$. Consider the family $\omega_s=\phi^*_s\omega$. This is a family of symplectic forms that all agree on $L$ and outside a neighborhood of $\Lambda_-$. Thus there is an isotopy $\Phi_s$ of $X$ that is fixed on $L$ and $\partial X$ such that $\Phi_s^*\omega_s=\omega_0=\omega$. Notice that $\Phi_s\circ\phi_s$ is a family of symplectomorphism of $(X,\omega)$ that covers the Lagrangian isotopy $L_s$. Moreover $(\Phi_s\circ \phi_s)_*v_-$ is a Liouville vector field that is tangent to $L$ near $Y_-$. 
\end{proof}
We are now ready for the main neighborhood theorem for a Lagrangian cobordism. This is a relative version of the standard Lagrangian neighborhood theorem. 
\begin{theorem}\label{relnbhd}
Let $L$ be a Lagrangian cobordism in $(X,\omega)$ from the Legendrian link $\Lambda_-$ in the contact manifold $\partial_- X$ to the Legendrian link $\Lambda_+$ in the contact manifold $\partial_+X$. Then there is a neighborhood $N$ of $L$ in $X$ and $N'$ of the zero section in $T^*L$ that are symplectomorphic. 
\end{theorem}
\begin{proof}
We first consider a boundary component $B$ of $L$. Suppose a neighborhood $N$ of $B$ is parameterized as $[0,1]\times S^1$ with $B=\{1\}\times S^1$. Use coordinates $(r,\theta)$ on this neighborhood. The cotangent space $T^*N$ is $[0,1]\times S^1\times \R^2$. If we use coordinates $(x,y)$ on $\R^2$, the Liouville form on $T^*N$ is given by $x\, dr + y\, d\theta$. Consider the vector field $v= r\partial_r + y\partial_y$. This vector field shows that $\{1\}\times S^1\times \R^2$ is a convex boundary of $T^*N$. Similarly, if $B$ is a negative boundary component of $L$ then $B\times \R^2$ will be a concave boundary component of $T^*N$.

A well-known result of Weinstein \cite[Proposition~4.2]{Weinstein91} says that there are neighborhoods of $\partial L$ in $X$ that are symplectomorphic to neighborhoods of $\partial L$ in $T^*L$ (by a symplectomorphism that takes the boundary components to the boundary components). As in the standard Lagrangian neighborhood theorem \cite[Section~3.4]{McDuffSalamon98}, we may extend this symplectomorphism to a diffeomorphism from a neighborhood of $L$ in $X$ to a neighborhood of the zero section in $T^*L$ so that $L$ is taken to the zero section and along the zero section $\omega$ is taken to the canonical symplectic structure on $T^*L$. Now a standard application of the Moser trick \cite[Section~3.2]{McDuffSalamon98} gives the claimed symplectomorphism. 
\end{proof}
With the above neighborhood theorem in place, we note that Lagrangian cobordisms can be isotoped to be symplectic. 
\begin{proposition}\label{perturbtosymp}
Let $L$ be a Lagrangian cobordism in $(X,\omega)$ from the Legendrian link $\Lambda_-$ in the contact manifold $\partial_- X$ to the Legendrian link $\Lambda_+$ in the contact manifold $\partial_+X$. If either $\Lambda_-$ or $\Lambda_+$ is non-empty, then $L$ may be perturbed to a symplectic cobordism $S$ from the transverse push-off of $\Lambda_-$ to the transverse push-off of $\Lambda_+$.
\end{proposition}
\begin{proof}
Using \cref{relnbhd} we only need to prove the zero-section in $T^*L$ can be so approximated. Given our hypothesis, there is a primitive $\beta$ for a volume form on $L$ such that $v_\beta$ points out of $L$ along $\Lambda_+$ and into $L$ along $\Lambda_-$, where $v_\beta$ is the vector field associated to $\beta$ by the symplectic form $d\beta$. Let $\lambda$ be the Liouville form on $T^*L$. Recall, this has the property that $\beta^*\lambda=\beta$ and therefore $\beta^*d\lambda=d\beta$. That is, the pull-back of the symplectic form on $T^*L$ by $\beta$ is a symplectic form on $L$. Thus the image of $\beta$ in $T^*L$ is a symplectic approximation of the zero-section by a symplectic submanifold and by scaling $\beta$ we may assume this approximation is arbitrarily close to the zero-section.

We are left to see that $S$ is a cobordism between the transverse push-offs of $\Lambda_-$ and $\Lambda_+$. To this end, notice that the vector field $v$ in the proof of \cref{relnbhd} induces a contact structure on $(T^*L)|_{\partial L}$ and in the coordinates given there we see the image of $\beta$ transverse to the contact structure and can be chosen to be $C^\infty$ close to $\Lambda_-\cup\Lambda_+$. 
\end{proof}

We are now ready to prove the symplectic analog of \cref{main}, that ribbon cobordisms can be approximated by symplectic cobordisms between transverse knots. 
\begin{thmMain3}
Given a ribbon cobordism $C$ in $[0,1]\times S^3$ from the link $K_-$ to the link $K_+$ with $K_-\neq\emptyset$, then there is some transverse realization $T_\pm$ of $K_\pm$ so that $C$ may be isotoped to a symplectic cobordism from $T_-$ to $T_+$. 
\end{thmMain3}

\begin{proof}[Proof of \cref{main3}]
Given a ribbon cobordism $C$ in $[0,1]\times S^3$ from the link $K_0$ to the link $K_1$, then we can use \cref{main} to find a Lagrangian cobordism $L$ between Legendrian realizations of the $K_i$. Now \cref{perturbtosymp} gives us the desired symplectic cobordism. 
\end{proof}

This allows us to prove the symplectic analog of \cref{main2}.
\begin{thmMain5}
Given a ribbon cobordism $C$ in $[0,1]\times S^3$ from the link $K_-$ to the link $K_+$ with $K_-\neq\emptyset$ and a transverse realization $T_+$ of $K_+$, then there is some transverse realization $T_-$ of $K_-$ so that $C$ may be isotoped to a symplectic cobordism from $T_-$ to $T_+$. 
\end{thmMain5}

\begin{proof}[Proof of \cref{main5}]
This proof is identical to the proof of \cref{main2} except we use \cref{main3} instead of \cref{main}.
\end{proof}

\section{Examples}\label{examples}

Given a ribbon cobordism from $K_-$ to $K_+$, \cref{main} tells us there is an isotopy of the cobordism to a Lagrangian cobordism from a Legendrian realization $\Lambda_-$ of $K_-$ to a Legendrian realization $\Lambda_+$ of $K_+$. Note that we do not know much about the Thurston-Bennequin numbers of $\Lambda_-$ and $\Lambda_+$. In this section we will explore some examples.

\begin{example}\label{ex:m61}
 \cref{fig:m61ribbon} gives a ribbon cobordism $C$ from the unknot to $m(6_1)$. Thus, \cref{main} tells us there exists a Legendrian unknot $\Lambda_-$ and Legendrian realization $\Lambda_+$ of $m(6_1)$ such that there is a decomposable Lagrangian cobordism $L$ from $\Lambda_-$ to $\Lambda_+$ which is smoothly isotopic to the ribbon cobordism. 
\begin{figure}[htbp]
\includegraphics{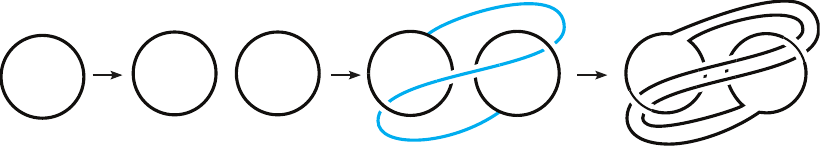}
\caption{A ribbon cobordism $C$ (in this case a filling) from the unknot to $m(6_1)$.}
\label{fig:m61ribbon}
\end{figure}
From \cite{Chantraine10}, we know $r(\Lambda_-)=r(\Lambda_+)$ and $\tbn(\Lambda_+)-\tbn(\Lambda_-)=-\chi(L)=0$, as $L$ is a filling of $\Lambda_+$. Thus $\tbn(\Lambda_+)=\tbn(\Lambda_-)$. Given that $\tbm(m(6_1))=-3$, we know that $\tbn(\Lambda_-)=\tbn(\Lambda_+)\leq-3$, which means that $\Lambda_-$ must be a Legendrian unknot that has been stabilized at least twice and that in this case it is impossible for $\Lambda_-$ to have maximal \tb number. Thus it is not always possible for the Lagrangian realization to have $\Lambda_-$ with maximal \tb number. However, in this case, it is possible for $\Lambda_+$ to have maximal \tb number. See \cref{fig:m61lagrangian} 
\begin{figure}[htbp]
\includegraphics{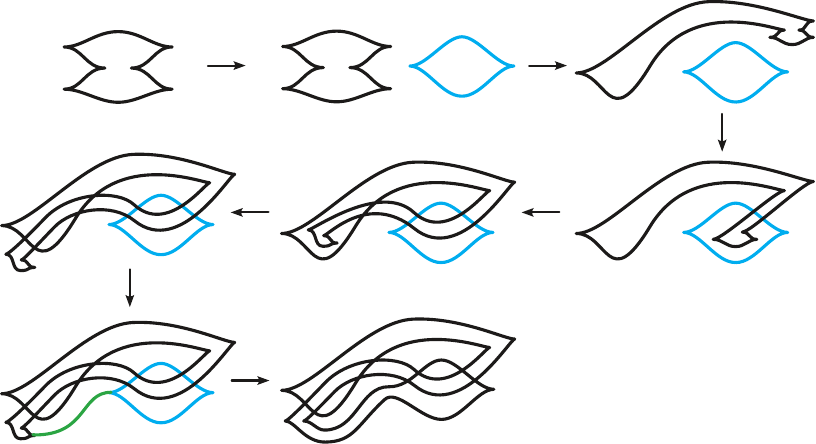}
\caption{A decomposable Lagrangian cobordism (in this case a filling) from an unknot to a max \tb $m(6_1)$ which is smoothly isotopic to the ribbon cobordism described in \cref{fig:m61ribbon}.}
\label{fig:m61lagrangian}
\end{figure}
for a description of a decomposable Lagrangian cobordism from a Legendrian unknot to a Legendrian $m(6_1)$ knot with $\tbn(\Lambda_+)=\tbm(m(6_1))=-3$ which is smoothly isotopic to the ribbon cobordism $C$.
\end{example}

The next example will show an infinite family of ribbon cobordism which have Lagrangian realizations with maximal \tb number for $\Lambda_+$. 

\begin{example}\label{ex:infiniteFamily}
Let $n\in\Z$ and define $P_n$ to be the pretzel knot $P(3,-3,n)$, see \cref{fig:infiniteFamily}. (In the case where $n=-1$, we see that $P_{-1}$ is $m(6_1)$ which was addressed in detail in \cref{ex:m61}.) Note that $P_n$ is a ribbon knot and so there exists a ribbon cobordism from the unknot to $P_n$. \cref{main} tells us there is a Legendrian representative $\Lambda_n$ of $P_n$ and a Legendrian unknot $U_n$ such that there is a decomposable Lagrangian cobordism from $U_n$ to $\Lambda_n$ which is smoothly isotopic to the original ribbon cobordism. However, \cref{main} does not tell us whether or not it is possible to find such a decomposable Lagrangian cobordism for a max \tb number representative of $P_n$. In \cite{Walker}, Walker shows that this is, in fact, possible for all $n$, where $\Lambda_n$ is depicted in \cref{fig:infiniteFamily}. Note that there are two places to perform a pinch move that will result in two Legendrian unknots, one with maximal \tb invariant. 

\begin{figure}[htbp]
\includegraphics{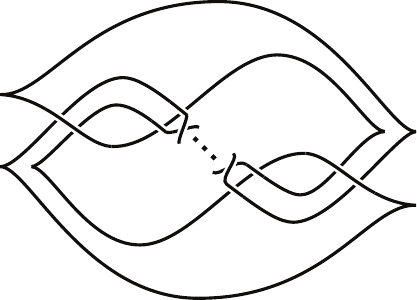}
\caption{Max Legendrian \tb number representatives of $P_n=P(3,-3,n)$ \cite{Walker}}
\label{fig:infiniteFamily}
\end{figure}
\end{example}

\begin{example}\label{ex:cableTorus}
Consider the two Legendrian $(2,3)$-cables of the trefoil knot in \cref{fig:cableTorus}. Both $\Lambda$ and $\Lambda'$ have rotation number 2 and Thurston-Bennequin number 5, while $\Lambda'$ destabilizes and $\Lambda$ does not. From \cite{Cables} we know there is only one Legendrian $(2,3)$-cable of the trefoil knot with rotation number $3$ and Thurston-Bennequin number $6$, so after positively stabilizing $\Lambda$ and $\Lambda'$ once the knots are Legendrian isotopic and so there exists a decomposable Lagrangian cobordism (in this case a concordance which is a trace of a Legendrian isotopy) from $S_+(\Lambda)$ to $S_+(\Lambda')$. However, from \cref{main4} we know there is a decomposable Lagrangian cobordism from a stabilization of $\Lambda$ to $\Lambda'$ (not stabilized). Using the LOSS invariant program \cite{LOSS} from Meyers, et. al., we see that $\hat\lambda^+(\Lambda)\neq0$ and $\hat\lambda^+(\Lambda')=0$. In \cite{BLW}, Baldwin, Lidman, and Wong showed that in such a case there is not a decomposable Lagrangian cobordism from $\Lambda$ to $\Lambda'$.

In summary, we know 
\begin{itemize}
\item there is a decomposable Lagrangian cobordism from $S_+(\Lambda)$ to $S_+(\Lambda')$,
\item there is not a decomposable Lagrangian cobordism from $\Lambda$ to $\Lambda'$,
\item but there is a decomposable Lagrangian cobordism from some stabilization of $\Lambda$ to $\Lambda'$ (not stabilized).
\end{itemize}
\begin{figure}[htbp]\tiny
    \begin{subfigure}{.45\textwidth}
    \begin{center}
    \includegraphics[width=.95\textwidth]{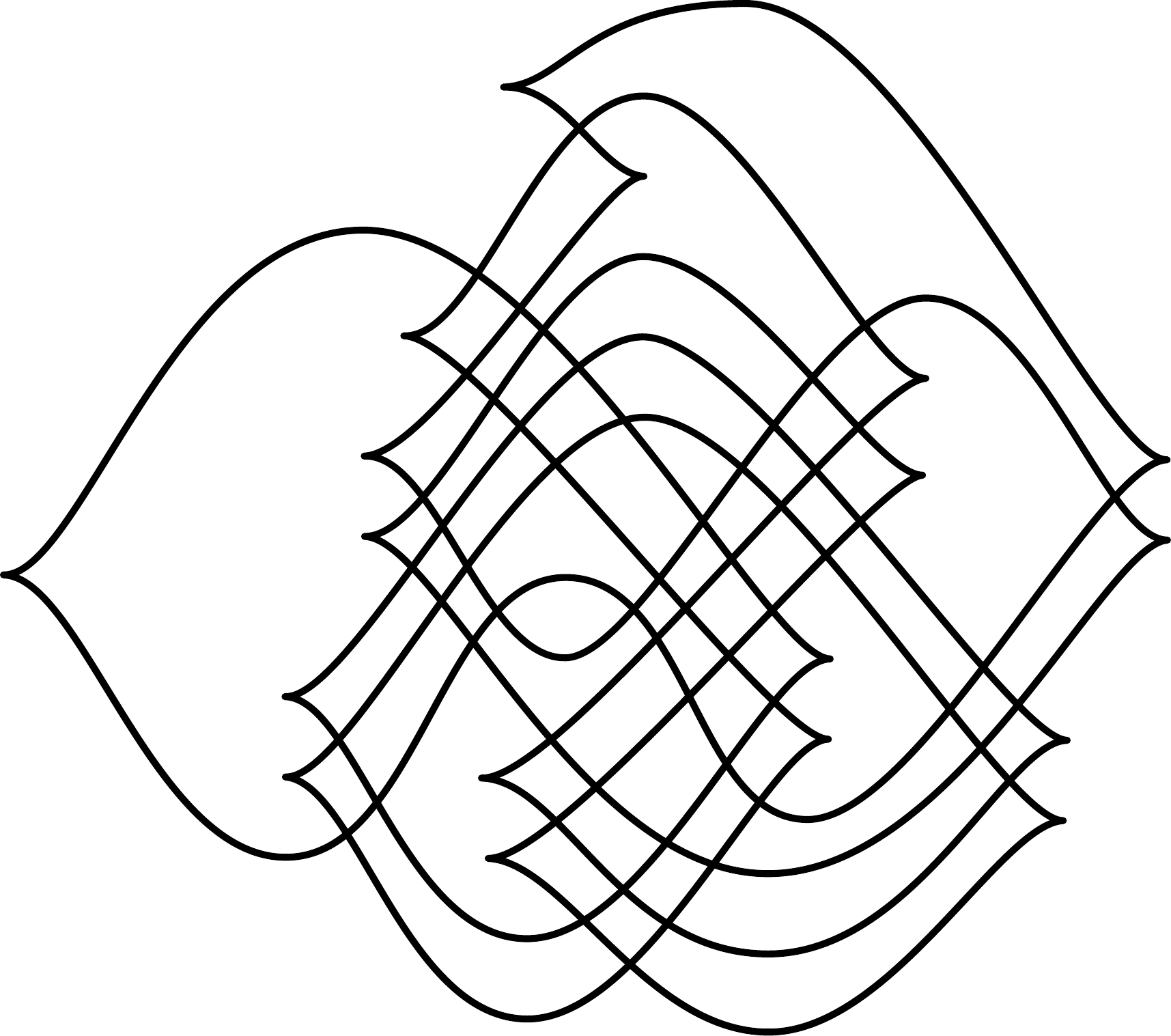}
    \caption{$\Lambda$}
    \end{center}
    \end{subfigure} 
    \begin{subfigure}{.45\textwidth}
    \begin{center}
    \includegraphics[width=.95\textwidth]{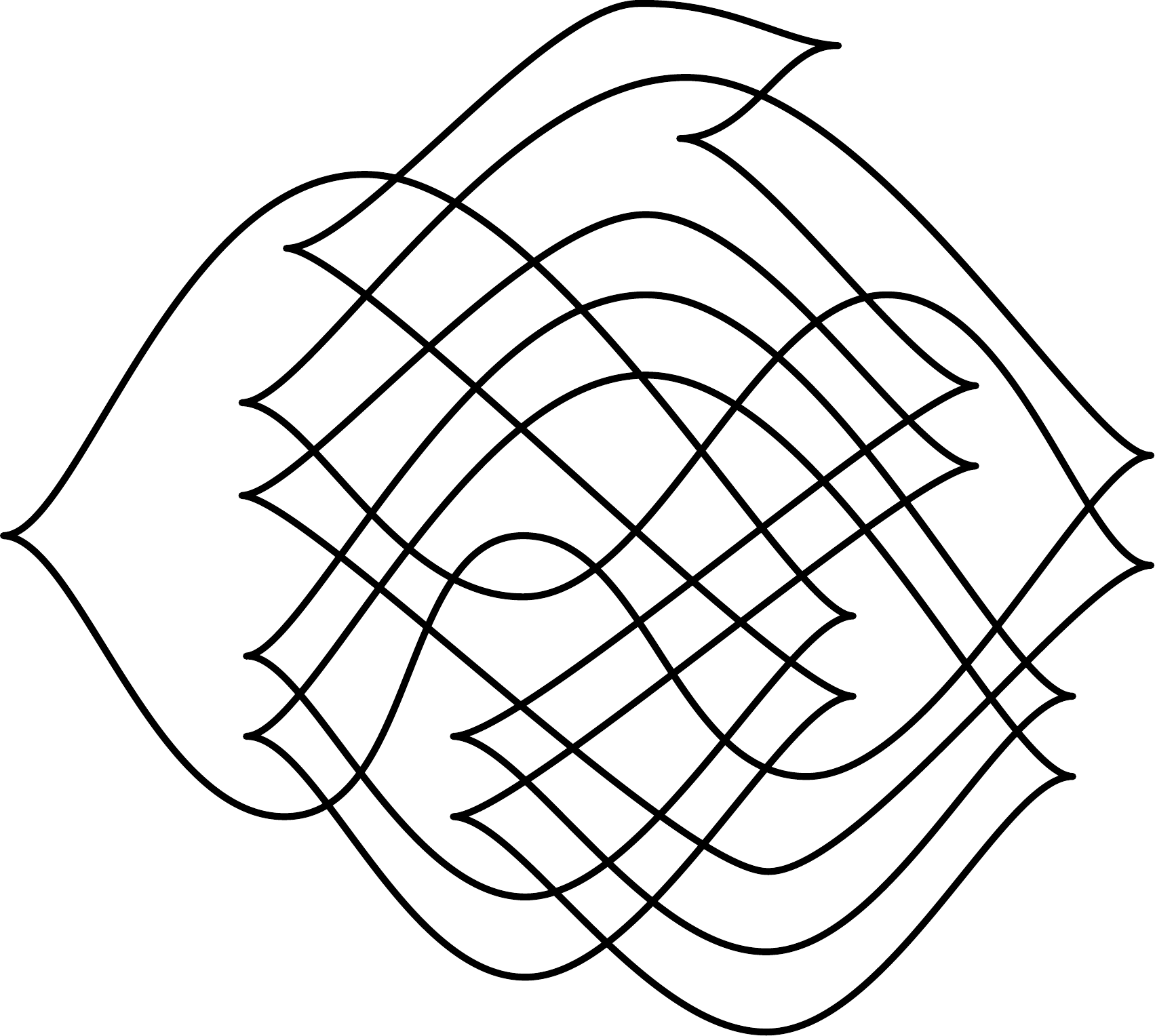}
    \caption{$\Lambda'$}
     \end{center}
    \end{subfigure}
    \caption{Two Legendrian $(2,3)$-cables of the $(2,3)$-torus knot.}
    \label{fig:cableTorus}
  \end{figure}

\end{example}

\bibliography{Bibliography}
\bibliographystyle{plain}

\end{document}